\documentclass[11pt,a4paper,fullpage]{article}
\usepackage[english]{babel}
\usepackage[utf8]{inputenc}
\usepackage{fancyhdr}
\usepackage[numbers]{natbib}

\usepackage{graphicx}
\usepackage{amsmath}
\usepackage{extsizes}
\usepackage{relsize}
\usepackage{multicol}
\usepackage{ragged2e}
\usepackage{fontenc}
\usepackage{mathcomp}
\usepackage{latexsym}
\usepackage{amssymb}
\usepackage{times}
\newtheorem{Teorema}{Theorem}[section]

\newtheorem{Definicija}[Teorema]{Definition}
\newtheorem{Posledica}[Teorema]{Corollary}
\newtheorem{Primer}[Teorema]{Example}
\newtheorem{Lema}[Teorema]{Lemma}
\newtheorem{Primedba}[Teorema]{Remark}
\usepackage{dsfont}
\numberwithin{equation}{section}
\hyphenchar\font=-1
\begin{document}
	\begin{center}
	{\LARGE Algebraic properties of the group inverse\par}	
\vspace{1cm}
	{\large Nikola Sarajlija}

\vspace{5mm}
	{\large \today\par}
\end{center}
\vspace{1cm}
\vspace{1cm}
\begin{abstract}
We investigate properties of the group inverse in rings with unit related to products and differences of idempotents, and thus we extend some results from \cite{DENG} to more general settings. We show that most part of \cite{DENG} is easily translated into the context of rings with unit, and that there are some difficulties in generalizing statements dealing with complementary idempotents.
\end{abstract}
\section{Introduction}
Field of generalized inverses has grown much in the last few decades. We shall study the group generalized inverse, which appears to have a central role in some applications \cite{MEYER, GOLUB, SOARES, MEYER2}, as well as in characterizations of some other generalized inverses \cite[section 3.3]{DJR}. Our main task will be to investigate which of some known results about group inverses for bounded linear operators still hold in a wider setting of rings with unit. 

Question of invertibility of $p-q$, $p,q$-idempotents, is especially important since it is connected with various other questions, such as integral equations, iterative methods in numerical linear algebra, signal processing and linear regression. For details, see introductory section in \cite{DENG} and the references mentioned there.

Throughout this text let $\mathcal{R}$ always stand for an arbitrary ring with unit. Its unit is denoted typically, by 1. 
For $a\in\mathcal{R}$, the group inverse of $a$ \cite{DJR} is every element $c\in\mathcal{R}$ that satisfies 
$$aca=a,\ cac=c,\ ac=ca.$$ 
Such  an element $c$ is unique in the case it exists, and we write $c=a^\#$. Moreover, $a^\#$ double commutes with $a$ and we say that $a$ is group invertible. This definition represents a natural extension of the definition of group inverse for square matrices, \cite{BENISRAEL}, \cite{CAMPBELL}, to arbitrary rings with unit. Group inverse can also be considered as a special case of the Drazin inverse. Namely, if Drazin inverse $a^D$ exists with $ind(a^D)\leq 1$, then $a$ is group invertible and $a^\#=a^D.$

In \cite{DENG} Deng has proved many interesting results concerning bounded linear idempotents $P, Q\in\mathcal{B}(\mathcal{H})$. Almost all of his proofs use the method of decomposing Hilbert space as $\mathcal{H}=\mathcal{R}(P)\oplus\mathcal{N}(P)$ for some idempotent $P$, and then studying the matrix decomposition of certain operators with respect to this decomposition. An appropriate analogue for this space decomposition in rings with unit is described as follows. If $a\in\mathcal{R}$ is group invertible, by $a^\pi=1-aa^\#$ we denote the so called spectral idempotent of $a$. It is easily checked that this idempotent generates the ring decomposition $\mathcal{R}=(a^\pi)^\circ\oplus a^\pi\mathcal{R}$, where $(a^\pi)^\circ$ and $a^\pi\mathcal{R}$ are the right annihilator and the right ideal of $a$, respectively. A quick reminder:
$$
x^\circ=\lbrace y\in\mathcal{R}:\ xy=0\rbrace,\ 
x\mathcal{R}=\lbrace xy:\ y\in\mathcal{R}\rbrace.
$$
Notice that with respect to this ring decomposition we can write $a$ in the unique way as $a=a_1\oplus 0$, where $a_1$ is invertible in the subring $(1-a^\pi)\mathcal{R}(1-a^\pi)$ with $a^\#$ being its inverse. Another way for writing this representation is $a=\bmatrix  a_1 & 0\\0 & 0\endbmatrix_{1-a^\pi}$, and this matrix representation will appear the most appropriate for our study. We introduce this form of representation in the subsequent section.
\section{Decompositions with respect to an idempotent}
Very important tool in work with generalized inverses is a representation of an arbitrary element $a\in\mathcal{R}$ in the matrix form using idempotent(s). There are two standard ways to represent an arbitrary element of a ring with unit in the matrix form: using one or two idempotents. Matrix representations based on the usage of two idempotents $p$ and $q$ are most often exploited in $C^*-$algebras and unital rings with involution, see for example novel work of Djordjevi\'{c} and Mihajlovi\'{c} \cite{MIHAJLOVIC}. Since our setting is much coarser than this, we choose using only one idempotent and describe this in the sequel. But first, let us remind ourselves how the usual representation using an idempotent looks like for the algebra of linear bounded operators.

If $\mathcal{H}$ is a Hilbert space, by $\mathcal{B}(\mathcal{H})$ we denote the algebra of linear bounded operators on $\mathcal{H}$. It is an algebra with a unit $I$, where $I$ is the identity operator. Let $P\in\mathcal{B}(\mathcal{H})$ be an idempotent, $P^2=P.$  Then the matrix representation of $P$ associated with the space decomposition $H=\mathcal{R}(P)\oplus\mathcal{N}(P)$ is 
\begin{equation}\label{Operatorska}
P=\begin{bmatrix}
I & 0 \\
0 & 0 
\end{bmatrix},
\end{equation} that is $P=I\oplus 0$. Arbitrary idempotent $Q\in\mathcal{B}(\mathcal{H})$, $Q^2=Q,$ has matrix representation \begin{equation}\label{Operatorska2}Q=\begin{bmatrix}
Q_1 & Q_2 \\
Q_3 & Q_4 
\end{bmatrix}
\end{equation} with respect to the space decomposition $\mathcal{H}=\mathcal{R}(P)\oplus\mathcal{N}(P).$ Since $Q^2=Q$ we obtain
$$\begin{bmatrix}
Q_1^2+Q_2Q_3 & Q_1Q_2+Q_2Q_4 \\
Q_3Q_1+Q_4Q_3 & Q_3Q_2+Q_4^2 
\end{bmatrix}=\begin{bmatrix}
Q_1 & Q_2 \\
Q_3 & Q_4 
\end{bmatrix}.$$

Analogously, let $p\in\mathcal{R}^\bullet$. It is easily checked that if $p$ is an idempotent, then $1-p$ is an idempotent as well. We call it complementary idempotent related to $p$ and sometimes denote it by $\overline{p}$. Let $a\in \mathcal{R}$. Then we write
$$
a=pap+pa(1-p)+(1-p)ap+(1-p)a(1-p)
$$
and use the notation
$$
a_{11}=pap,\quad a_{12}=pa(1-p),\quad a_{21}=(1-p)ap,\quad a_{22}=(1-p)a(1-p) .
$$
Hence, every projection $p\in\mathcal{R}$ induces the representation of an arbitrary element
$a\in\mathcal{R}$ given by the following matrix
$$ a=\bmatrix  pap & pa(1-p)\\(1-p)ap &
(1-p)a(1-p)\endbmatrix_p = \bmatrix  a_{11}(p)&
a_{12}(p)\\a_{21}(p)&a_{22}(p)\endbmatrix_p.
$$
As a very important feature, we must say that multiplication in $\mathcal{R}$ agrees with multiplication of matrices. That is, if $a$ and $b$ are represented by matrices $\bmatrix  a_{11}&
a_{12}\\a_{21}&a_{22}\endbmatrix_p$ and $\bmatrix  b_{11}&
b_{12}\\b_{21}&b_{22}\endbmatrix_p$, than product $ab$ can be calculated by the usual matrix multiplication, namely $$ab=\bmatrix a_{11}b_{11}+a_{12}b_{21}&a_{11}b_{12}+a_{12}b_{22}\\a_{21}b_{11}+a_{22}b_{21}& a_{21}b_{12}+a_{22}b_{22}\endbmatrix_p.$$
\noindent If $p$ is well-known, then we use $a_{ij}$ instead of $a_{ij}(p),\
i,j\in\{1,2\}.$

If $q\in\mathcal{R}^\bullet$ is an arbitrary idempotent, then we obviously have
\begin{equation}\label{glavna}
p=\begin{bmatrix}p&0\\0&0\end{bmatrix}_p,\quad q=\begin{bmatrix}q_{1}&q_{2}\\q_{3}&q_{4}\end{bmatrix}_p,\quad p,q_1,q_2,q_3,q_4\in\mathcal{R},
\end{equation}
and this is the desired representation of an element in $\mathcal{R}$ that is an appropriate analogue of (\ref{Operatorska}),  (\ref{Operatorska2}). By analogy with the similar situation for operators, from $q^2=q$ we get \begin{equation}\label{glavna2}\begin{bmatrix}
q_1^2+q_2q_3 & q_1q_2+q_2q_4 \\
q_3q_1+q_4q_3 & q_3q_2+q_4^2 
\end{bmatrix}_p=\begin{bmatrix}
q_1 & q_2 \\
q_3 & q_4 
\end{bmatrix}_p,
\end{equation} where $q$ is given in (\ref{glavna}).

We will give a very important observation in a connection with representation (\ref{glavna}). Namely, if we compare $P$ from (\ref{Operatorska}) with $p$ from $(\ref{glavna}$), then we see an important difference. In (\ref{Operatorska}) we have an instance of the identity operator, which is the unit of the algebra $\mathcal{B}(H)$. However, in (\ref{glavna}) instead of $1$, the unit of the ring $\mathcal{R}$, we have $p$ at the corresponding position. This will have a huge impact on our work in the following way: many analogues of theorems appearing in \cite{DENG} will have proofs that no longer hold unless some additional assumptions are added. For that reason, we will modify those assumptions so that they exclude some of the idempotents $q$ represented as in (\ref{glavna}) (this is different from \cite{DENG} where we had corresponding results holding \textit{for all} idempotents $Q$ of an appropriate representation). Since $p$ is the unit of the subring $p\mathcal{R}p$, the easiest way to accomplish this is to leisurely assume that idempotents $q\in p\mathcal{R}p$. In other words, instead of (\ref{glavna}), we would suggest
\begin{equation}\label{glavna3}
p=\begin{bmatrix}p&0\\0&0\end{bmatrix}_p,\quad q=\begin{bmatrix}q_{1}&q_{2}\\q_{3}&q_{4}\end{bmatrix}_p,\quad p,q_1,q_2,q_3,q_4\in p\mathcal{R}p.
\end{equation}
However, this solution appears to be too strict in a sense that assuming it all results from the following two sections trivially hold. So, in the sequel we relax this condition $q\in p\mathcal{R}p$ in an appropriate way (Theorems \ref{teorema2.1}, \ref{teorema2.2} and their corollaries...)

As a supplement, let us just say that this matrix representation of elements in a ring is useful in many other areas of this theory, and not only for group inverses. For example, the generalized Drazin inverse can be finely represented using this tool. It is well-known \cite{DJR} that $a\in \mathcal{R}^d$ (the set of all elements being generalized Drazin invertible) can be represented in the
following matrix form
$$
a=\bmatrix a_1&0\\0&a_2\endbmatrix_p
$$
with respect to $p=aa^d=1-a^\pi$, where $a^\pi$ is the spectral idempotent of $a$,
$a_1$ is invertible in the subring $p\mathcal{R} p$ and $a_2$ is quasinilpotent in the subring
$(1-p)\mathcal{R} (1-p)$. Notice that $a_2$ is quasinilpotent in the ring $\mathcal{R}$ as well. Then, generalized Drazin
inverse of $a$ is given by
$$
a^d=\bmatrix  [a_1]_{p\mathcal{R}p}^{-1}&0\\0&0\endbmatrix_p.
$$
\section{Basic results}
First we state the following result:
\begin{Lema}\label{lema2.1} Every idempotent $a\in\mathcal{R}$ is group invertible with $a^{\#}$=$a$. If $a$ and $b$ are group invertible and commute, then
$$(ab)^\#=b^\#a^\#=a^\#b^\#,\ a^\#b=ba^\#\  and\  b^\#a=ab^\#.$$
If $ab$ and $ba$ are group invertible, then\\
\begin{equation}\label{Klajn}(ab)^\#=a[(ba)^\#]^2b.\end{equation}
\end{Lema}
Proof of this lemma is familiar in the literature and it follows from the fact $ab=ba$, so we omit it. See \cite{DENG} for the operator case.
\begin{Lema}
\label{lema2.2}
(\cite{OLESKY}, Theorem 2.2 for the matrix case, and \cite{DENG}, Lemma 2.2 for the operator case]) Let $m\in\mathcal{R}$ have the matrix form
$$m=\begin{bmatrix}
0 & a \\
b & 0 
\end{bmatrix}_p.$$
Then $m$ is group invertible if and only if $ab$ and $ba$ are group invertible and
$$(ab)^\pi a=0,\quad b(ab)^\pi=0.$$
In this case,
$$m^\#=\begin{bmatrix}
0 & (ab)^\#a \\
b(ab)^\# & 0 
\end{bmatrix}_p=\begin{bmatrix}
0 & a(ba)^\# \\
(ba)^\#b & 0 
\end{bmatrix}_p.$$
\end{Lema}
\textbf{Proof. }Assume that $ab$ and $ba$ are group invertible and $(ab)^\pi a=0, b(ab)^\pi=0$. Let $x=\begin{bmatrix}
0 & (ab)^\#a \\
b(ab)^\# & 0 
\end{bmatrix}_p.$ From $b(ab)^\#a=(ba)^\#ba$ we get
$$mx=xm=(ab)^\#ab\oplus b(ab)^\#a=(ab)^\#ab\oplus(ba)^\#ba.$$
So,
$$xmx=\begin{bmatrix}
0 & (ab)^\#a \\
b(ab)^\# & 0 
\end{bmatrix}_p\begin{bmatrix}
(ab)^\#ab & 0 \\
0 & b(ab)^\#a 
\end{bmatrix}_p=x.$$
From $(ab)^\pi a=0$ and $b(ab)^\pi=0$ we get
$$mxm=\begin{bmatrix}
(ab)^\#ab & 0 \\
0 & b(ab)^\#a 
\end{bmatrix}_p\begin{bmatrix}
0 & a \\
b & 0 
\end{bmatrix}_p=\begin{bmatrix}
0 & a \\
b & 0 
\end{bmatrix}_p=m.$$
Hence, $m$ is group invertible and $m^\#=x$. 

Assume that $m$ is group invertible. Then $m^2=\begin{bmatrix}
ab & 0 \\
0 & ba 
\end{bmatrix}_p$ is group invertible and, therefore, $ab$ and $ba$ are group invertible. Moreover, group invertibility of $m$ implies that $m\mathcal{R}=m^2\mathcal{R}$, that is $ ab\mathcal{R}=a\mathcal{R}$ and $ba\mathcal{R}=b\mathcal{R}$. The equality $ab(ab)^\#ab=ab$ implies that $ab(ab)^\#x=x$ for all vectors $x\in ab\mathcal{R}$. Thus $ab(ab)^\#x=x$ for all $x\in a\mathcal{R}$, and therefore $ab(ab)^\#a=a$.\\
Similarly, we obtain the second desired equality.
$\square$

Let $p$ and $q$ be two idempotents. Next, we provide some equivalent conditions for the existence of group inverses of $p - q$, $\overline{p} - q$ and $pq - qp$ under the assumption that $p$ and $q$ are given by (\ref{glavna}). As already said, we demand that $p$ and $q$ satisfy some extra conditions (whose analogues do not exist in \cite{DENG}) that we shall call \textit{additional assumptions}. 
\begin{Teorema}\label{teorema2.1}\textit{Let $p$ and $q$ be idempotents given by (\ref{glavna}).}\\[1mm]
(i) Let us assume that $pq_1=q_1p=q_1,\ pq_2=q_2$ and $q_3p=q_3$ (additional assumptions). Then $p - q$ is group invertible if and only if $p - q_1$ and $q_4$ are group 
invertible and
\begin{equation}\label{sixth}
\begin{split}
\left(p - q_1\right)^\# q_2 = q_2 q_4^\#, & \quad q_2 q_4^{\pi} = \left(p - q_1\right)^{\pi} q_2 = 0,\\
    q_3\left(p - q_1\right)^\# = q_4^\# q_3, & \quad  q_4^{\pi} q_3 = q_3\left(p - q_1\right)^{\pi} = 0.
\end{split}
\end{equation}
In this case,
\begin{equation}\label{seventh}
    \left(p - q\right)^\# = \begin{bmatrix}
    \left(p - q_1\right)^\#\left(p - q_1\right) & - \left(p - q_1\right)^\# q_2 \\ - q_4^\# q_3 & - q_4^\# q_4 
    \end{bmatrix}_p;
\end{equation}
(ii) Let us assume that $\overline{p}q_4=q_4\overline{p}=q_4,\ q_2\overline{p}=q_2$ and  $\overline{p}q_3=q_3$ (additional assumptions). Then $\overline{p} - q$ is group invertible if and only if $q_1$ and $\overline{p} - q_4$ are group 
invertible and
\begin{equation}\label{eight}
\begin{split}
\left(\overline{p} - q_4\right)^\# q_3 = q_3 q_1^\#, & \quad q_3 q_1^{\pi} = \left(\overline{p} - q_4\right)^{\pi} q_3 = 0,\\
    q_2\left(\overline{p} - q_4\right)^\# = q_1^\# q_2, &\quad  q_1^{\pi} q_2 = q_2\left(\overline{p} - q_4\right)^{\pi} = 0.
\end{split}
\end{equation}
In this case,
\begin{equation}
    \left(\overline{p} - q\right)^\# = \begin{bmatrix}
    - q_1^\# q_1 & - q_1^\# q_2 \\
    - \left(\overline{p} - q_4\right)^\# q_3 &  \left(\overline{p} - q_4\right)^\#  \left(\overline{p} - q_4\right)
    \end{bmatrix}_p;
\end{equation}
(iii) Let us assume that additional assumptions from $(i)$ and $(ii)$ hold:
\begin{equation*}
\begin{split}
pq_1=q_1p=q_1,\ & pq_2=q_2\  and\  q_3p=q_3,\\
\overline{p}q_4=q_4\overline{p}=q_4,\ & q_2\overline{p}=q_2\  and\  \overline{p}q_3=q_3.
\end{split}
\end{equation*}
If $p - q$ and $\overline{p} - q$ are group invertible, then $pq - qp$ is group invertible and 
\begin{equation}\label{ten}
    \left(pq - qp\right)^\# = \begin{bmatrix}
    0 & - q_1^\# \left(p - q_1\right)^\# q_2  \\
    q_3 q_1^\# \left(p - q_1\right)^\# & 0
    \end{bmatrix}_p.
\end{equation}
\end{Teorema}
\textbf{Proof.} (i) Assume that $p - q_1$, $q_4$ are group invertible and expressions in (\ref{sixth}) hold. Let us denote by $x$ the right side of (\ref{seventh}). Then
$$
    \left(p - q\right)x =  \begin{bmatrix}
    p - q_1 & - q_2 \\
    - q_3 & - q_4 
    \end{bmatrix}_p \begin{bmatrix}
    \left(p - q_1\right)^\#\left(p - q_1\right) & - \left(p - q_1\right)^\# q_2 \\ - q_4^\# q_3 & - q_4^\# q_4 
    \end{bmatrix}_p
$$
$$
    =  \begin{bmatrix}
    \left(p - q_1\right)\left(p - q_1\right)^\# & 0 \\ 0 &  q_4 q_4^\# 
    \end{bmatrix}_p 
$$
because $\left(p - q_1\right)\left(p - q_1\right)^\# q_2 = q_2 = q_2 q_4^\# q_4$ and $q_3\left(p - q_1\right)\left(p - q_1\right)^\# = q_3 = q_4 q_4^\# q_3$ in view of (\ref{sixth}),
$$\begin{aligned}p - q_1 + q_2 q_4^\# q_3 = p - q_1 + \left(p - q_1\right)^\# q_2 q_3 =\left(p - q_1\right)\left(p - q_1\right)^\#\left(p - q_1\right) +\\ \left(p - q_1\right)^\#\left(p - q_1\right) q_1 = \left(p - q_1\right)^\#\left(p - q_1\right),\end{aligned}$$
and $q_3\left(p - q_1\right)^\# q_2 + q_4 = q_3 q_2 q_4^\# + q_4 = \left(p - q_4\right)q_4 q_4^\# + q_4 = q_4 q_4^\#$ in view of (\ref{sixth}) and (\ref{glavna2}). We have also used the additional assumptions. We will not stress use of additional assumptions anymore. In a similar way we get
$$
    x\left(p - q\right) = \begin{bmatrix}
    \left(p - q_1\right)^\#\left(p - q_1\right) & - \left(p - q_1\right)^\# q_2 \\ - q_4^\# q_3 & - q_4^\# q_4 
    \end{bmatrix}_p \begin{bmatrix}
    p - q_1 & - q_2 \\
    - q_3 & - q_4 
    \end{bmatrix}_p 
$$
$$
    =  \begin{bmatrix}
    \left(p - q_1\right)\left(p - q_1\right)^\# & 0 \\ 0 &  q_4 q_4^\# 
    \end{bmatrix}_p 
$$
because $ \left(p - q_1\right)^\# \left(p - q_1\right) q_2 = q_2 = q_2 q_4^\# q_4 =  \left(p - q_1\right)^\# q_2 q_4$ and \\
$q_4^\# q_3 \left(p - q_1\right) = q_3 \left(p - q_1\right)^\# \left(p - q_1\right) = q_3 = q_4 q_4^\# q_3$ in view of (\ref{sixth}), $p - q_1 +  \left(p - q_1\right)^\# q_2 q_3 = p - q_1 +  \left(p - q_1\right)^\# \left(p - q_1\right) q_1 =  \left(p - q_1\right)^\# \left(p - q_1\right)$ and $q_4^\# q_3 q_2 + q_4 = q_4^\# q_4 \left(p - q_1\right) + q_4 = q_4 q_4^\#$ in view of (\ref{sixth}) and (\ref{glavna2}). So
\begin{equation}\label{twelve}
    x \left(p - q\right) = \left(p - q\right) x = \left(p - q_1\right)\left(p - q_1\right)^\# \oplus q_4 q_4^\#.
\end{equation}
Now,
$$
    \left(p - q\right) x\left(p - q\right) = \begin{bmatrix}
    p - q_1 & - q_2 \\
    - q_3 & - q_4 
    \end{bmatrix}_p \begin{bmatrix}
    \left(p - q_1\right)\left(p - q_1\right)^\# & 0 \\ 0 &  q_4 q_4^\# 
    \end{bmatrix}_p 
$$
$$    
    = \begin{bmatrix}
    p - q_1 & - q_2 \\
    - q_3 & - q_4 
    \end{bmatrix}_p 
$$
and
$$
    x \left(p - q\right) x = \begin{bmatrix}
    \left(p - q_1\right)^\#\left(p - q_1\right) & - \left(p - q_1\right)^\# q_2 \\ - q_4^\# q_3 & - q_4^\# q_4 
    \end{bmatrix}_p \begin{bmatrix}
    \left(p - q_1\right)\left(p - q_1\right)^\# & 0 \\ 0 &  q_4 q_4^\# 
    \end{bmatrix}_p 
$$
$$
    = \begin{bmatrix}
    \left(p - q_1\right)^\#\left(p - q_1\right) & - \left(p - q_1\right)^\# q_2 \\ - q_4^\# q_3 & - q_4^\# q_4 
    \end{bmatrix}_p .
$$
Hence, $p - q$ is group invertible and $\left(p - q\right)^\#$ has representation (\ref{seventh}).\\[1mm]
(ii) Observe that $\overline{p} - q = 1 - p - q = - \left[p - \left(1 - q\right)\right]$. Thus group invertibility of $\overline{p} - q$ is equivalent to group invertibility of $p - \left(1 - q\right)$. If $p$ and $q$ are represented as in (\ref{glavna}), then $p - q$ has an obvious representation by means of (\ref{glavna}), and we can apply item (i) in an evident manner.\\[1mm]
(iii) If $p - q$ and $\overline{p} - q$ are group invertible, by (\ref{glavna2}) and Lemma \ref{lema2.1},
$$
    \left(q_2 q_3\right)^\# = q_1^\#\left(p - q_1\right)^\# = \left(p - q_1\right)^\# q_1^\#, \quad  
$$
$$    
   \left(q_3 q_2\right)^\# = q_4^\#\left(p - q_4\right)^\# = \left(p - q_4\right)^\# q_4^\#.
$$
By (\ref{sixth}) and (\ref{eight}),
$$
     \left(q_2 q_3\right) \left(q_2 q_3\right)^\#  q_2 = \left(p - q_1\right)\left(p - q_1\right)^\# q_1 q_1^\# q_2 = \left(p - q_1\right)\left(p - q_1\right)^\# q_2 = q_2
$$
and
$$
    \left(q_3 q_2\right) \left(q_3 q_2\right)^\#  q_3 = \left(\overline{p} - q_4\right)\left(\overline{p} - q_4\right)^\# q_4 q_4^\# q_3 = \left(\overline{p} - q_4\right)\left(\overline{p} - q_4\right)^\# q_3 = q_3.
$$
Hence, $q_2 q_3$, $q_3 q_2$ are group invertible, $\left(q_2 q_3\right)^{\pi} q_2 = 0 $ and $\left(q_3 q_2\right)^{\pi} q_3 = 0$. By Lemma \ref{lema2.2}, $pq - qp = \begin{bmatrix}
    0 & q_2 \\
    - q_3 & 0
    \end{bmatrix}_p$ is group invertible and (\ref{ten}) follows directly by Lemma \ref{lema2.2}. 

Assume that $p - q$ is group invertible. Then $\left(p - q\right)^2$ is group invertible. Since $\left(p - q\right)^2 = p + q - pq - qp$, by using representation (\ref{glavna}) we get $\left(p - q\right)^2 = \left(p - q_1\right) \oplus q_4$. Hence $p - q_1$, $q_4$ are group invertible and
$$
    \left(p - q\right)^\# = \left[ \left(p - q\right)^2\right]^\# \left(p - q\right) = \left(p - q\right) \left[ \left(p - q\right)^2\right]^\#
$$
$$
     = \begin{bmatrix}
    p - q_1 & 0 \\
    0 & q_4 
    \end{bmatrix}_p^\# 
      \begin{bmatrix}
    p - q_1 & - q_2 \\
    - q_3 & - q_4 
    \end{bmatrix}_p = 
       \begin{bmatrix}
    p - q_1 & - q_2 \\
    - q_3 & - q_4 
    \end{bmatrix}_p
    \begin{bmatrix}
    p - q_1 & 0 \\
    0 & q_4 
    \end{bmatrix}_p^\# 
$$
$$
    = \begin{bmatrix}
    \left(p - q_1\right)^\#\left(p - q_1\right) & - \left(p - q_1\right)^\# q_2 \\ - q_4^\# q_3 & - q_4^\# q_4 
    \end{bmatrix}_p = \begin{bmatrix}
    \left(p - q_1\right)\left(p - q_1\right)^\# & - q_2 q_4^\# \\ - q_3 \left(p - q_1\right)^\# & - q_4 q_4^\# 
    \end{bmatrix}_p.
$$
Comparing two sides of the above equation, we have 
\begin{equation}\label{eleventh}
    \left(p - q_1\right)^\# q_2 = q_2 q_4^\# \quad and \quad q_3\left(p - q_1\right)^\# = q_4^\# q_3.
\end{equation}
From $\left(p - q\right)\left(p - q\right)^\#\left(p - q\right) = \left(p - q\right)$ and 
$$
    \left(p - q\right)\left(p - q\right)^\# = \left(p - q\right)^2\left[ \left(p - q\right)^2\right]^\# = \left(p - q_1\right)\left(p - q_1\right)^\# \oplus q_4 q_4^\#,
$$
we get
$$
     \begin{bmatrix}
    \left(p - q_1\right)\left(p - q_1\right)^\# & 0 \\ 0 &  q_4 q_4^\# 
    \end{bmatrix}_p 
    \begin{bmatrix}
    p - q_1 & - q_2 \\
    - q_3 & - q_4 
    \end{bmatrix}_p = \begin{bmatrix}
    p - q_1 & - q_2 \\
    - q_3 & - q_4 
    \end{bmatrix}_p.
$$
Hence,
$$
    q_2 = \left(p - q_1\right)\left(p - q_1\right)^\# q_2, \quad\quad q_3 = q_4 q_4^\# q_3.
$$
The first equality of (\ref{eleventh}) yields $\left(p - q_1\right)\left(p - q_1\right)^\# q_2 = \left(p - q_1\right) q_2 q_4^\#$. Moreover, (\ref{glavna2}) leads to $\left(p - q_1\right) q_2 q_4^\# = q_2 q_4 q_4^\#$. The second equality of (\ref{eleventh}) yields $q_4 q_4^\# q_3 = q_4 q_3\left(p - q_1\right)^\#$. Moreover, (\ref{glavna2}) leads to \\$q_4 q_3\left(p - q_1\right)^\# = \left(q_3 - q_3 q_1\right)\left(p - q_1\right)^\# = q_3\left(p - q_1\right)\left(p - q_1\right)^\#$. $\square$

\begin{Posledica}\label{posledica2.1} \textit{Let $p$ and $q$ be idempotents given by (\ref{glavna})}.\\[1mm]
(i) If $p - q$ is group invertible with $$pq_1=q_1p=q_1,pq_2=q_2\  and\  q_3p=q_3,$$ then $\left(p - q\right)^{\pi} = \left(p - q_1\right)^{\pi} \oplus q_4^{\pi}$;\\[1mm]
(ii) If $\overline{p} - q$ is group invertible with $$\overline{p}q_4=q_4\overline{p}=q_4,\ q_2\overline{p}=q_2\  and\   \overline{p}q_3=q_3,$$ then $\left(\overline{p} - q\right)^{\pi} = q_1^{\pi} \oplus \left(\overline{p} - q_4\right)^{\pi}$;\\[1mm]
(iii) \textit{If $p - q$ and $\overline{p} - q$ are group invertible with 
\begin{equation*}\begin{split}pq_1=q_1p=q_1,pq_2=q_2&\  and\  q_3p=q_3,\\ 
\overline{p}q_4=q_4\overline{p}=q_4,\  q_2\overline{p}=q_2&\  and\  \overline{p}q_3=q_3,
\end{split}
\end{equation*} then $\left(pq - qp\right)^{\pi} = \left(q_1^2 - q_1\right)^{\pi} \oplus \left(q_4 - q_4^2\right)^{\pi}$}.
\end{Posledica}

\noindent\textbf{Proof. } Use (\ref{glavna2}) and Theorem \ref{teorema2.1}. $\square$

Next topic we are interested in is the product of idempotents. We state the second main theorem of this article, again, of course, under some additional assumptions.

\begin{Teorema}\label{teorema2.2} \textit{Let $p$ and $q$ be idempotents given by (\ref{glavna})}. Suppose that in (i-ii) $p$ commutes with $q_1, q_2$, and in (iii-iv) $p$ commutes with $q_3,q_4$. Then:\\[1mm]
(i) $pq$ is group invertible if and only if $pq_1$ is group invertible and $(pq_1)^{\pi}pq_2 = 0$;\\[1mm]
(ii) $\overline{p}q$ is group invertible if and only if $pq_4$ is group invertible and $(pq_4)^{\pi}pq_3 = 0$;\\[1mm]
(iii) $p\overline{q}$ is group invertible if and only if $p-pq_1$ is group invertible and $(p-pq_1)^{\pi}pq_2 = 0$;\\[1mm]
(iv) $\overline{p}\hspace{1mm}\overline{q}$ is group invertible if and only if $p - pq_4$ is group invertible and \\$\left(p - pq_4\right)^{\pi}pq_3 = 0$.
\end{Teorema}
\textbf{Proof.} We only prove item (i), items (ii)-(iv) can be proved in the same way. Let $p$ and $q$ be given by (\ref{glavna}).
Then
$$
    pq = \begin{bmatrix}
    pq_1 & pq_2 \\
    0 & 0
    \end{bmatrix}_p, \quad \left(pq\right)^2 = \begin{bmatrix}
    pq_1^2 & pq_1q_2 \\
    0 & 0
    \end{bmatrix}_p, \quad and \quad \left(pq\right)^3 = \begin{bmatrix}
    pq_1^3 & pq_1^2q_2 \\
    0 & 0
    \end{bmatrix}_p,
$$
with respect to the ring decomposition $\mathcal{R} = p\mathcal{R} \oplus p^\circ$ and with the help of additional assumptions. Since $pq$ is group invertible, ind$(pq) = 1$. So
$$
    pq\mathcal{R} = (pq)^2\mathcal{R} =  (pq)^3\mathcal{R}=\cdots\quad and \quad  (pq)^\circ = ((pq)^2) ^\circ= ((pq)^3)^\circ=\cdots
$$
For every $x \in(pq_1^2)^\circ, x \oplus 0 \in((pq)^2) ^\circ= (pq)^\circ$. So $pq_1 x = 0$ and $(pq_1^2)^\circ \subset (pq_1)^\circ$. Since $(pq_1)^\circ \subset (pq_1^2)^\circ$ is trivial, we get $(pq_1^2)^\circ = (pq_1)^\circ$. For every $z \in pq_1\mathcal{R}$ there exists $x \in p\mathcal{R}$ such that 
$$
    \begin{bmatrix}
    pq_1 &p q_2 \\
    0 & 0
    \end{bmatrix}_p \begin{bmatrix}
    x \\
    0 
    \end{bmatrix}_p = \begin{bmatrix}
    z \\
    0 
    \end{bmatrix}_p \in pq\mathcal{R} = (pq)^3\mathcal{R}.
$$
So, there exist  $u \in p\mathcal{R}$ and $v \in p^\circ$ such that 
$$
    \begin{bmatrix}
    pq_1^3 & pq_1^2q_2 \\
    0 & 0
    \end{bmatrix}_p \begin{bmatrix}
    u \\
    v 
    \end{bmatrix}_p = \begin{bmatrix}
    z \\
    0 
    \end{bmatrix}_p.
$$
Therefore, $z = pq_1^3u +pq_1^2 q_2 v = pq_1^2\left(q_1 u + q_2 v\right) \in pq_1^2\mathcal{R}$ and $pq_1\mathcal{R} \subset pq_1^2\mathcal{R}$. Since $pq_1^2\mathcal{R} \subset pq_1\mathcal{R} $ is trivial, we get $pq_1^2\mathcal{R} = pq_1\mathcal{R}$. Hence, ind$(pq_1) \leq 1$, i.e. $(pq_1)^\#$ exists. Now, let
\begin{equation}\label{thirteen}
    x = \begin{bmatrix}
    pq_1^\# & \left((pq_1)^\#\right)^2pq_2 \\
    0 & 0
    \end{bmatrix}_p.
\end{equation}
It is trivial to check that $pqx = xpq$, $xpqx = x$ and $(pq)^3x = (pq)^2$. These imply that $x$ is the Drazin inverse of $pq$ and ind($pq) \leq 2$. Since $pq$ is group invertible, $(pq)^2x = pq$. We get
$$
    \begin{bmatrix}
    pq_1^2 &p q_1q_2 \\
    0 & 0
    \end{bmatrix}_p \begin{bmatrix}
    pq_1^\# & \left((pq_1)^\#\right)^2pq_2 \\
    0 & 0
    \end{bmatrix}_p = \begin{bmatrix}
    pq_1 & pq_2 \\
    0 & 0
    \end{bmatrix}_p,
$$
which implies that $pq_1 (pq_1)^\# pq_2 = pq_2$, i.e. $(pq_1)^{\pi} pq_2 = 0$.

On the other hand, if $pq_1$ is group invertible and $(pq_1)^{\pi}pq_2 = 0$, it is easy to check that (\ref{thirteen}) is group inverse of $pq$. $\square$

Here is another theorem which is trivial for the $q$ as in (\ref{glavna3}), and holds for the $q$ as in (\ref{glavna}) under additional assumptions.
\begin{Teorema}\textit{Let $p$ and $q$ be idempotents given by (\ref{glavna}). Suppose that $p$ commutes with $p_i$,\ $i=1,2,3$. Then $pq - qp$ is group invertible if and only if $pq_i\left(1 - q_i\right)p, i = 1,4$ are group invertible, $[pq_1\left(1 - q_1\right)p]^{\pi} pq_2 = 0$ and $pq_3[pq_4\left(1 - q_4\right)p]^{\pi}= 0$.}
\end{Teorema}
\textbf{Proof.} By
(\ref{glavna}), $pq - qp = \begin{bmatrix}
    0 & q_2p \\
    - q_3p & 0
    \end{bmatrix}_p$. By (\ref{glavna2}), $pq_2 q_3p = pq_1\left(1 - q_1\right)p$ and $pq_3 q_2p = pq_4\left(1 - q_4\right)p$. The result follows immediately by Lemma \ref{lema2.2}. $\square$

\begin{Primedba}
Note that
$$
    pq + qp = - (p + q) (\overline{p} - q) = - (\overline{p} - q)(p + q).
$$
By Lemma \ref{lema2.1}, if $p + q$ and $\overline{p} - q$ are group invertible, then $pq + qp$ is group invertible. $\triangle$
\end{Primedba}
\section{Further results}
In this section, we will obtain explicit representations for group inverses of products and differences of projections. First, we give the following definition.
\begin{Definicija}\label{definicija3.1}Let $p$ and $q$ be idempotents such that $p - q$ is group invertible. Define elements $f,g$ and $h$ as
\begin{equation}\label{fourteen}
    f = p\left(p - q\right)^\#, \quad g = \left(p - q\right)^\#p, \quad h = \left(p - q\right)^\#\left(p - q\right).
\end{equation}
\end{Definicija}

In what follows, we will have many statements concerning not only idempotents $p$ and $q$, but also their complementary idempotents $\overline{p}$ and $\overline{q}$ as well. For the algebra of linear bounded operators $\mathcal{B}(H)$, representations of $P$ and $\overline{P}$ used in \cite{DENG} are the following:
\begin{equation*}
P=\begin{bmatrix}I & 0\\0 & 0\end{bmatrix},\quad\overline{P}=\begin{bmatrix}0 & 0\\0 & I\end{bmatrix}.\\
\end{equation*}
In the ring $\mathcal{R}$, however, representations of $p$ and $\overline{p}$ are:
\begin{equation*}
p=\begin{bmatrix}p & 0\\0 & 0\end{bmatrix},\quad\overline{p}=\begin{bmatrix}0 & 0\\0 & \overline{p}\end{bmatrix}.\\
\end{equation*}
Notice the difference between $\overline{P}$ and $\overline{p}$. In the bottom right corner of the matrix representation, $\overline{P}$ has operator $I$, and $\overline{p}$ has element $\overline{p}$ at the same position. $I$ is an identity in the algebra $\mathcal{B}(H)$, but $\overline{p}$ is not an identity in the ring $\mathcal{R}$ nor in the ring $p\mathcal{R}p$. This will leave some significant consequences in the structure of this section. Namely, not all of the statements in \cite{DENG} related to complementary idempotents will have their analogues in the setting $\mathcal{R}$. This is the reason why some of the mentioned analogues will be left out, and we will not emphasize this fact. For example, in part (i) of the next theorem one would want to have $f=(p-q)^\#\overline{q}$ (the analogue of \cite[Theorem 3.1 (i)]{DENG}). However, this equality does not hold even under additional assumptions of Theorem \ref{teorema3.1}. If we add extra conditions $q_2p=q_4^\#q_4p=0$ then we would get this equality for $f$. However, for $g=\overline{q}(p-q)^\#$ \cite[Theorem 3.1 (ii)]{DENG} we would need some extra conditions different from the previous ones and so on... This is the reason for us to go no further with additional conditions than the ones stated in Theorem \ref{teorema3.1}. However, in part (ii) of Theorem \ref{teorema3.3} we will be forced to make an exception.
\begin{Teorema}\label{teorema3.1} \textit{Let $p$ and $q$ be idempotents given by (\ref{glavna}) such that $p - q$ is group invertible with $pq_1=q_1p=q_1,\ pq_2=q_2,\ q_3p=p$, and let $f,g$ and $h$ be given by Definition \ref{definicija3.1}. Then $f,g$ and $h$ are idempotents and}\\[1mm]
(i) $ f\mathcal{R} = ph\mathcal{R} = (p-q_1)\mathcal{R}$.\\[1mm]
(ii) $g^\circ = (ph)^\circ = p^\circ \oplus (p - q_1)^\circ$.
\end{Teorema}
\textbf{Proof.} From Theorem \ref{teorema2.1}(i) it follows
\begin{equation}\label{fifteen}
    f = \begin{bmatrix}
    \left(p - q_1\right)^\#\left(p - q_1\right) & - \left(p - q_1\right)^\# q_2 \\ 0 & 0 
    \end{bmatrix}_p, \quad g = \begin{bmatrix}
    \left(p - q_1\right)^\#\left(p - q_1\right) & 0 \\ - q_3\left(p - q_1\right)^\# & 0 
    \end{bmatrix}_p.
\end{equation}
and
\begin{equation}\label{sixteen}
    h = \left(p - q_1\right)^\#\left(p - q_1\right) \oplus q_4^\# q_4.
\end{equation}
So, $f,g$ and $h$ are idempotents,
$$
    f\mathcal{R} = \left[\left(p - q_1\right)^\#\left(p - q_1\right)\right]\mathcal{R} = \left[p - q_1\right]\mathcal{R} = ph\mathcal{R},
$$
and
$$
    g^\circ= \left[\left(p - q_1\right)^\#\left(p - q_1\right)\right]^\circ \oplus p^\circ = \left(p - q_1\right)^\circ \oplus p^\circ= (ph)^\circ.
$$
$\square$

Let $\mathcal{S}$ and $\mathcal{T}$ be subrings of $\mathcal{R}$ such that $\mathcal{R}=\mathcal{S}\oplus\mathcal{T}.$ If there exists idempotent $p\in\mathcal{R}$ with $p\mathcal{R} = \mathcal{S},\ p^\circ = \mathcal{T}$, we denote it by $p_{\mathcal{S},\mathcal{T}}$.
\begin{Posledica}\label{posledica3.1} Let $p$ and $q$ be given by (\ref{glavna}) with $pq_1=q_1p=q_1,\ pq_2=q_2,\ q_3p=p$, and let $f,g$ and $h$ be given by Definition \ref{definicija3.1}. If $p - q$ is invertible, then\\[1mm]
(i) $f = p \left(p - q\right)^{-1} = p_{p\mathcal{R},q\mathcal{R}}$;\\[1mm]
(ii) $g = \left(p - q\right)^{-1}p = p_{q^\circ,p^\circ}$;\\[1mm]
(iii) $\overline{p}f = 0$;\\[1mm]
(iv) $g \overline{p} = 0$.
\end{Posledica}
\textbf{Proof.} It is enough to note that if $p - q$ is invertible we will have that $p - q_1$, $q_4$ are invertible and $h = 1$. $\square$
\begin{Posledica}\label{posledica3.2} Let $p$ and $q$ be idempotents given by (\ref{glavna}) such that $p - q$ is group invertible with $pq_1=q_1p=q_1,\ pq_2=q_2,\ q_3p=p$, and let $f$ and $g$ be given by Definition \ref{definicija3.1}. Then
\begin{equation}\label{seventeen}
    fg =  \left(p - q_1\right)^\# \oplus 0. 
\end{equation}
\end{Posledica}
\textbf{Proof.} Combine (\ref{glavna2}), (\ref{fifteen})-(\ref{sixteen}). $\square$
Further algebraic manipulations yield some extra properties of $f,g$ and $h$.
\begin{Teorema} Let $p$ and $q$ be idempotents given by (\ref{glavna}) such that $p - q$ is group invertible with $pq_1=q_1p=q_1,\ pq_2=q_2,\ q_3p=p$, and let $f,g$ and $h$ be given by Definition \ref{definicija3.1}. Then\\[1mm]
(i) $fp = pg = ph = hp$;\\[1mm]
(ii) $qhq = qh = hq = hqh$.\\[1mm]
\end{Teorema}
\textbf{Proof.} (i) By (\ref{glavna}), (\ref{fifteen}) and (\ref{sixteen}), we get that
$$
    fp = pg = ph = hp = \left(p - q_1\right)\left(p - q_1\right)^\# \oplus 0.
$$
(ii) If $p - q$ is group invertible, by Theorem \ref{teorema2.1}(i), $q_2q_4^{\pi} = \left(p - q_1\right)^{\pi} q_2 = 0$ and $q_4^{\pi} q_3 = q_3\left(p - q_1\right)^{\pi} = 0$. Then
$$
    qh = \begin{bmatrix}
    q_1 & q_2 \\ q_3 &  q_4 
    \end{bmatrix}_p \begin{bmatrix}
    \left(p - q_1\right)\left(p - q_1\right)^\# & 0 \\ 0 &  q_4 q_4^\# 
    \end{bmatrix}_p = \begin{bmatrix}
    \left(q_1 - q_1^2\right)\left(p - q_1\right)^\# & q_2 \\ q_3 &  q_4  
    \end{bmatrix}_p
$$
and
$$
     hq = \begin{bmatrix}
    \left(p - q_1\right)\left(p - q_1\right)^\# & 0 \\ 0 &  q_4 q_4^\# 
    \end{bmatrix}_p \begin{bmatrix}
    q_1 & q_2 \\ q_3 &  q_4 
    \end{bmatrix}_p = \begin{bmatrix}
    \left(q_1 - q_1^2\right)\left(p - q_1\right)^\# & q_2 \\ q_3 &  q_4  
    \end{bmatrix}_p.
$$
So $qh = hq $ and
$$
    hqh = \begin{bmatrix}
    \left(p - q_1\right)\left(p - q_1\right)^\# & 0 \\ 0 &  q_4 q_4^\# 
    \end{bmatrix}_p\begin{bmatrix}
    \left(q_1 - q_1^2\right)\left(p - q_1\right)^\# & q_2 \\ q_3 &  q_4  
    \end{bmatrix}_p = qh.
$$
Moreover, from
$$
    q_1 \left[\left(q_1 - q_1^2\right)\left(p - q_1\right)^\#\right] + q_2q_3 = q_1^2 \left[\left(p - q_1\right)\left(p - q_1\right)^\#\right] + q_1\left(p - q_1\right) \quad by\ (\ref{glavna2})
$$
$$
    = q_1^2 \left[\left(p - q_1\right)\left(p - q_1\right)^\#\right] + q_1\left(p - q_1\right) \left[\left(p - q_1\right)\left(p - q_1\right)^\#\right] 
$$
$$    
    = \left(q_1 - q_1^2\right)\left(p - q_1\right)^\#
$$
and
$$
    q_3 \left[\left(q_1 - q_1^2\right)\left(p - q_1\right)^\#\right] + q_4q_3 = q_3q_1 \left[\left(p - q_1\right)\left(p - q_1\right)^\#\right] + q_3\left(p - q_1\right) \quad by\ (\ref{glavna2})
$$
$$
    = q_3q_1 \left[\left(p - q_1\right)\left(p - q_1\right)^\#\right] + q_3\left(p - q_1\right) \left[\left(p - q_1\right)\left(p - q_1\right)^\#\right] 
$$
$$
    = q_3 \left[\left(p - q_1\right)\left(p - q_1\right)^\#\right] \quad by\ (\ref{sixth})
$$
$$
    = q_3.
$$
we have
$$
    qhq =  \begin{bmatrix}
    q_1 & q_2 \\ q_3 &  q_4 
    \end{bmatrix}_p  \begin{bmatrix}
    \left(q_1 - q_1^2\right)\left(p - q_1\right)^\# & q_2 \\ q_3 &  q_4  
    \end{bmatrix}_p 
$$
$$    
    =  \begin{bmatrix}
    \left(q_1 - q_1^2\right)\left(p - q_1\right)^\# & q_2 \\ q_3 &  q_4  
    \end{bmatrix}_p = hq.
$$
$\square$

As a crown of preceding considerations we get the following result.

\begin{Teorema}\label{teorema3.3} Let $p$ and $q$ be idempotents given by (\ref{glavna}) such that $p - q$ is group invertible with $pq_1=q_1p=q_1,\ pq_2=q_2,\ q_3p=p$, and let $f,g$ and $h$ be given by Definition \ref{definicija3.1}. Then\\[1mm]
(i) $\left(p - q\right)^\# = f+g-h;$\\[1mm]
(ii) If we demand $q_2p=q_4p=0$, then $$\left(p + q\right)^\# = \left(p - q\right)^\#\left(p + q\right)\left(p - q\right)^\# = \left(2g - h\right)\left(f + g - h\right)$$ if and only if $ph = p$.
\end{Teorema}
\textbf{Proof.} Item (i) follows directly by Theorem \ref{teorema2.1}(i) and relations in (\ref{fifteen})-(\ref{sixteen}). So we only need to prove the item (ii). Here we have not only the additional assumptions ($pq_1=q_1p=q_1,\ pq_2=q_2,\ q_3p=q_3$), but also extra conditions $q_2p=0,\ q_4p=0$. These extra conditions ensure that $(p-q)^\#(\overline{p}+\overline{q})=(p+q)(p-q)^\#$, and in this way we provide that the proof of this theorem follows the same path as the one presented in \cite[Theorem 3.3 (ii)]{DENG}. For convenience of a reader we pursue with full proof. 
Denote by $x = \left(p - q\right)^\#\left(p + q\right)\left(p - q\right)^\#$. By the previous argument we have
$$
    \left(p + q\right)x = \left(p + q\right)\left(p - q\right)^\#\left(p + q\right)\left(p - q\right)^\#
$$
$$
    = \left(p - q\right)^\#\left(\overline{p} + \overline{q}\right)\left(p + q\right)\left(p - q\right)^\#
$$
$$
    = \left(p - q\right)^\#\left(p - q\right)^2\left(p - q\right)^\#
$$
$$
    = \left(p - q\right)\left(p - q\right)^\# = h,
$$
$$
    x \left(p + q\right) = \left(p - q\right)^\#\left(p + q\right)\left(p - q\right)^\#\left(p + q\right)
$$
$$
    = \left(p - q\right)^\#\left(p + q\right)\left(\overline{p} + \overline{q}\right)\left(p - q\right)^\#
$$
$$
    = \left(p - q\right)^\#\left(p - q\right)^2\left(p - q\right)^\#
$$
$$
    = \left(p - q\right)\left(p - q\right)^\# = h,
$$
and 
$$
    x \left(p + q\right)x = x\left(p - q\right)\left(p - q\right)^\#
$$
$$
    = \left(p - q\right)^\#\left(p + q\right)\left(p - q\right)^\#\left(p - q\right)\left(p - q\right)^\#
$$
$$
   = \left(p - q\right)^\#\left(p + q\right)\left(p - q\right)^\# = x.
$$
Hence, $\left(p + q\right)^\# = x$ if and only if $\left(p + q\right)x\left(p + q\right) = p + q$. Next, we prove that $\left(p + q\right)x\left(p + q\right) = p + q$ if and only if $ph = p$. We divide the proof in two steps:

First, if $ph = p$, by (\ref{sixth}) and (\ref{sixteen}), we have that $p - q_1$ is invertible, $qh = q$ and 
$$
    \left(p + q\right)x\left(p + q\right) = \left(p + q\right)h = p + q.
$$
Second, note that
$$
    \left(p + q\right)x\left(p + q\right) = \left(p + q\right)h
$$
$$
    = \begin{bmatrix}
    p + q_1 & q_2 \\ q_3 &  q_4  
    \end{bmatrix}_p\begin{bmatrix}
    \left(p - q_1\right)^\#\left(p - q_1\right) & 0 \\ 0 &  q_4^\# q_4 
    \end{bmatrix}_p
$$
$$
    = \begin{bmatrix}
    \left(p + q_1\right)\left(p - q_1\right)^\#\left(p - q_1\right) & q_2 \\ q_3 &  q_4  
    \end{bmatrix}_p. \quad by\  (\ref{sixth})
$$
If $\left(p + q\right)x\left(p + q\right) = p + q$, then $\left(p+ q_1\right)\left(p - q_1\right)^\#\left(p - q_1\right) = p+ q_1$. Since $p - q_1$ as an element acting on $p\mathcal{R}$ is group invertible,
$$
    p\mathcal{R} = \left(\left(p - q_1\right)^{\pi}\right)^\circ
 \oplus \left(\left(p - q_1\right)^{\pi}\right)\mathcal{R} = \left(p - q_1\right)\mathcal{R} \oplus\left(p - q_1\right)^\circ.
$$
For every $x \in \left(p - q_1\right)^\circ$, we have $q_1 x = x$ and $2x = \left(p + q_1\right)x=$\\$ \left(p + q_1\right)\left(p- q_1\right)^\#\left(p - q_1\right)x = 0.$

It follows that $\left(p - q_1\right)^\circ = \{0\}$ and so $(p-q_1)\mathcal{R}=\mathcal{R}$ and $p - q_1$ is invertible. Hence $ph = p$.

Now, by the definition of group inverse, we know $$x = \left(p - q\right)^\#\left(p + q\right)\left(p - q\right)^\# = \left(p + q\right)^\#.$$ Moreover, by item (i) of this theorem,
$$
    \left(p + q\right)^\# = \left(p - q\right)^\#\left(p + q\right)\left(p - q\right)^\#
$$
$$
    = \left(p - q\right)^\# p \left(p - q\right)^\# + \left(p - q\right)^\# q \left(p - q\right)^\#
$$
$$
    = g\left(f + g - h\right) + \left(g - h\right)\left(f + g - h\right)
$$
$$
    = \left(2g - h\right)\left(f + g - h\right).
$$
$\square$
\begin{Primedba}
Condition $ph = p$ in Theorem \ref{teorema3.3}, item (ii) is necessary to assure $\left(p + q\right)x\left(p + q\right) = p + q$.
This can be demonstrated in the following example.
\end{Primedba}
\begin{Primer} Define idempotents $p,q$ on $\oplus_{i = 1}^4 \mathds{R}$ ($\mathds{R}$ is the set of real numbers) by
$$
    p = 1 \oplus 1 \oplus 0 \oplus 0 \quad and \quad q = 1 \oplus 0 \oplus 1 \oplus 0.
$$
Then
$$
    \left(p - q\right)^\# = 0 \oplus 1 \oplus - 1 \oplus 0,
$$
$$
    \left(p + q\right)^\# = \frac{1}{2} \oplus 1 \oplus 1 \oplus 0,
$$
$$
    h = \left(p - q\right)\left(p - q\right)^\# = 0 \oplus 1 \oplus 1 \oplus 0,
$$
$$
    x = \left(p - q\right)^\#\left(p + q\right)\left(p - q\right)^\# = 0 \oplus 1 \oplus 1 \oplus 0.
$$
In this case, $ph \neq p, \left(p + q\right)x\left(p + q\right) \neq p + q$ and, hence $\left(p + q\right)^\# \neq x$. $\triangle$
\end{Primer}
\begin{Posledica} \label{posledica3.4}\cite[Theorem 2.2]{RAKOCEVIC}. \textit{Let $p$ and $q$ be idempotents given by (\ref{glavna}) such that $pq_1=q_1p=q_1,\ pq_2=q_2,\ q_3p=q_3$, $f = p_{p\mathcal{R},q\mathcal{R}}$ and $g = p_{q^\circ ,p^\circ}$. If $p - q$ is invertible with $q_2p=q_4p=0$ in addition, then}\\[1mm]
(i) $\left(p + q\right)^{-1} = \left(p - q\right)^{-1}\left(p + q\right)\left(p - q\right)^{-1}$;\\[1mm]
(ii) $\left(p - q\right)^{-1} = \left(p + q\right)^{-1}\left(p - q\right)\left(p + q\right)^{-1}$;\\[1mm]
(iii) $\left(p - q\right)^{-1} = f+g - 1$;\\[1mm]
(iv) $\left(p + q\right)^{-1} = \left(2g - 1\right)\left(f + g - 1\right)$.
\end{Posledica}

If $p - q$ is invertible, Harte gave that \cite[Theorem 7.5.1]{HARTE}
$$
    \left(1 - pqp\right)^{-1} = 1 - p + p\left(p - q\right)^{-2}.
$$
In fact, we can obtain some similar results on group inverses.
\begin{Teorema}\label{teorema3.4}  Let $p$ and $q$ be idempotents given by (\ref{glavna}) such that $p - q$ is group invertible with $pq_1=q_1p=q_1,\ pq_2=q_2,\ q_3p=p$, and let $f,g$ and $h$ be given by Definition \ref{definicija3.1}. Then\\[1mm]
(i) $\left(p - pqp\right)^\# = fg$;\\[1mm]
(ii) $\left(p - pq\right)^\# = [fg]^2\overline{q}$;\\[1mm]
(iii) $\left(p - qp\right)^\# = \overline{q}[fg]^2$.
\end{Teorema}
\textbf{Proof.} (i) From (\ref{glavna}), $p - pqp = \left(p - q_1\right) \oplus 0$, hence by Corollary \ref{posledica3.2},
$$
    \left(p - pqp\right)^\# = \left(p - q_1\right)^\# \oplus 0 = fg.
$$
(ii) Note that
$$
    \left(p - pq\right)^\# = \begin{bmatrix}
    \left(p - q_1\right)^\# & -\left[\left(p - q_1\right)^\# \right]^2 q_2 \\ 0 &  0  
    \end{bmatrix}_p
$$
$$
    = \begin{bmatrix}
    \left[\left(p - q_1\right)^\# \right]^2 & 0 \\ 0 &  0  
    \end{bmatrix}_p \begin{bmatrix}
    p - q_1 & -q_2 \\ -q_3 &  \overline{p} - q_4  
    \end{bmatrix}_p 
$$
$$
    = \left(fg\right)^2 \overline{q}.
$$
(iii) Note that 
$$
    \left(p - qp\right)^\# = \begin{bmatrix}
    \left(p - q_1\right)^\# & 0 \\ -q_3\left[\left(p - q_1\right)^\# \right]^2 &  0  
    \end{bmatrix}_p
$$
$$
    = \begin{bmatrix}
    p - q_1 & -q_2 \\ -q_3 &  \overline{p} - q_4  
    \end{bmatrix}_p  \begin{bmatrix}
    \left[\left(p - q_1\right)^\# \right]^2 & 0 \\ 0 &  0  
    \end{bmatrix}_p 
$$
$$
    = \overline{q} \left(fg\right)^2 .\quad \square
$$

\noindent\textbf{Acknowledgments}\\[3mm]
\hspace*{6mm}This work was supported by the Ministry of Education, Science and Technological Development of the Republic of Serbia under Grant No. 451-03-9/2021-14/ 200125\\[3mm]
\textbf{Declaration of interest}\\[3mm]
\hspace*{6mm}I wish to confirm that there are no known conflicts of interest associated with this publication and
there has been no significant financial support for this work that could have influenced its outcome other than listed above. I confirm that there are
no other persons who satisfied the criteria for authorship but are not listed.
I confirm that there are no impediments to publication with respect to intellectual property.\\[3mm]
\noindent\textbf{Data availability statement}\\[3mm]
\hspace*{6mm}The datasets generated during and/or analysed during the current study are available from the corresponding author on reasonable request.

\begin{flushleft}
Corresponding author:

Nikola Sarajlija, \\
nikola.sarajlija@dmi.uns.ac.rs\\
+381637461610\\
University of Novi Sad, Faculty of Sciences, Serbia
\end{flushleft}
\end{document}